\documentclass[12pt,twoside]{article}

\usepackage{amssymb}
\usepackage{amsmath}
\usepackage{bbm}
\usepackage{mathrsfs}
\usepackage{xypic}

\makeatletter
\gdef\th@mychange{\normalfont\slshape
   \def\@begintheorem##1##2{\item
        [\hskip\labelsep \theorem@headerfont ##2. ##1  \,--\!--\!--\!--  ]}%
 \def\@opargbegintheorem##1##2##3{%
   \item[\hskip\labelsep \theorem@headerfont ##2. ##1\ {\upshape(}##3{\upshape)}. \,-----  ]}}
\makeatother

\RequirePackage{theorem}
\theoremstyle{mychange}

{\theorembodyfont{\rmfamily}\newtheorem{ttt}{}[section]}
{\theorembodyfont{\rmfamily}\newtheorem{nota}[ttt]{Notation.}}
{\theorembodyfont{\rmfamily}\newtheorem{remark}[ttt]{Remark.}}
{\theorembodyfont{\rmfamily}\newtheorem{rems}[ttt]{Remarks.}}
{\theorembodyfont{\rmfamily}\newtheorem{ex}[ttt]{Example.}}
{\theorembodyfont{\rmfamily}}

{\theorembodyfont{\itshape}\newtheorem{subl}[ttt]{Sublemma.}}
{\theorembodyfont{\itshape}\newtheorem{lem}[ttt]{Lemma.}}
{\theorembodyfont{\itshape}\newtheorem{coro}[ttt]{Corollary.}}
{\theorembodyfont{\itshape}\newtheorem{prop}[ttt]{Proposition.}}
{\theorembodyfont{\itshape}\newtheorem{theo}[ttt]{Theorem.}}

{\theorembodyfont{\rmfamily}\newtheorem{meto}[ttt]{Approach}}
{\theorembodyfont{\rmfamily}\newtheorem{imp}[ttt]{\em The improvement.}}

\newcounter{abc}
\newenvironment{abc}{\begin{list}{\rm \alph{abc}) }{\usecounter{abc} \leftmargin=0.0pt \labelsep=0.0pt \listparindent=0.0pt \labelwidth=0.0pt \parsep=\smallskipamount \itemsep=0.0pt \topsep=0.0pt \partopsep=\smallskipamount}}{\end{list}}
\newcounter{iii}
\newenvironment{iii}{\begin{list}{\rm \roman{iii}) }{\usecounter{iii} \leftmargin=0.0pt \labelsep=0.0pt \listparindent=0.0pt \labelwidth=0.0pt \parsep=\smallskipamount \itemsep=0.0pt \topsep=0.0pt \partopsep=\smallskipamount}}{\end{list}}

\newcommand{\Gal}{\mathop{\text{\rm Gal}}\nolimits}
\newcommand{\Spec}{\mathop{\text{\rm Spec}}\nolimits}
\newcommand{\Pic}{\mathop{\text{\rm Pic}}\nolimits}
\newcommand{\Picb}{\mathop{\text{\bf Pic}}\nolimits}
\newcommand{\coker}{\mathop{\text{\rm coker}}\nolimits}
\renewcommand{\sp}{\mathop{\text{\rm sp}}\nolimits}

\newcommand{\im}{\mathop{\text{\rm im}}\nolimits}
\newcommand{\Frob}{\mathop{\text{\rm Frob}}\nolimits}
\newcommand{\rk}{\mathop{\text{\rm rk}}\nolimits}

\newcommand{\limdir}{\mathop{\underrightarrow{\rm lim}}\limits}

\newcommand{\Pb}{{\bf P}}

\newcommand{\et}{\text{\'et}}
\newcommand{\sing}{\text{sing}}

\newcommand{\calI}{\mathscr{I}}
\newcommand{\calK}{\mathscr{K}}
\newcommand{\calL}{\mathscr{L}}
\newcommand{\calO}{\mathscr{O}}

\newcommand{\bbC}{{\mathbbm C}}
\newcommand{\bbF}{{\mathbbm F}}
\newcommand{\bbN}{{\mathbbm N}}
\newcommand{\bbQ}{{\mathbbm Q}}
\newcommand{\bbZ}{{\mathbbm Z}}

\newcommand{\frakm}{{\mathfrak m}}

\newcommand{\nr}{{\rm nr}}

\newcommand{\br}{ }
\newcommand{\brr}{, }

\def\rightend#1#2{{%
 \leavevmode\nobreak\hskip .5em plus 1fil
 \penalty600 \hskip 0pt plus -1filll
 \vadjust{}\nobreak\hskip 0pt plus 1filll%
 #1\parfillskip=#2\relax \par}}
\def\eop{\ifmmode\rule[-22pt]{0pt}{1pt}\ifinner\tag*{$\square$}\else\eqno{\square}\fi\else\rightend{$\square$}{0pt}\fi}

\author{Andreas-Stephan Elsenhans and J\"org Jahnel}
\date{}

\title{The Picard group of a $K3$~surface and \\its reduction modulo~$p$}

\begin{document}
\renewcommand{\thefootnote}{\fnsymbol{footnote}}

\maketitle

\begin{abstract}
We~present a method to compute the geometric Picard rank of a
$K3$~surface
over~$\bbQ$.
Contrary~to a widely held belief, we show it is possible to verify Picard
rank~$1$
using reduction only at a single~prime. Our~method is based on deformation theory for invertible~sheaves.
\end{abstract}

\section{Introduction}

\begin{ttt}
For~$K3$~surfaces,
the Picard group is a highly interesting~invariant. In~general, it is
isomorphic
to~$\bbZ^n$
for
some~$n = 1, \ldots, 20$.
A~generic
$K3$~surface
over~$\bbC$
has Picard
rank~$1$.
Nevertheless,~the first explicit examples of
$K3$~surfaces
over~$\bbQ$
with geometric
Picard~rank~$1$
were constructed by R.~van Luijk~\cite{vL} as late as~2004. Van~Luijk's method is based on reduction
modulo~$p$.
It~works as~follows.
\end{ttt}

\begin{meto}[{\rm van Luijk}{}]
\label{vLm}
Let~$S$
be
a~$K3$~surface.

\begin{iii}
\item
At a
place~$p$
of good reduction, the Picard
group~$\Pic(S_{\overline{\bbQ}})$
of the surface injects into the
Picard
group~$\smash{\Pic(S_{\overline{\bbF}_{\!p}})}$
of its reduction
modulo~$p$.
\item
On~its part,
$\smash{\Pic(S_{\overline{\bbF}_{\!p}})}$
injects into the second \'etale cohomology group
$\smash{H_\et^2(S_{\overline{\bbF}_{\!p}}\!, \bbQ_l(1))}$.\vspace{0.1mm}
\item
Only roots of unity can arise as eigenvalues of the
Frobenius on the image
of~$\Pic(S_{\overline{\bbF}_{\!p}})$
in~$\smash{H_\et^2(S_{\overline{\bbF}_{\!p}}\!, \bbQ_l(1))}$.
The~number of eigenvalues of this form is therefore an upper bound for the Picard rank
of~$S_{\overline{\bbF}_{\!p}}$.
One~may compute the eigenvalues
of~$\Frob$
by counting the points
on~$S$,
defined
over~$\bbF_{\!p}$
and some finite~extensions.
\end{iii}\smallskip

\noindent
Doing~this for one prime, one obtains an upper bound
for~$\rk \Pic(S_{\overline{\bbF}_{\!p}})$
which is always~even. The~Tate conjecture asserts that this bound is actually~sharp. For~proving that the Picard rank
over~$\overline{\bbQ}$
is equal
to~$1$,
the best that could happen is to find a prime which yields an upper bound
of~$2$.\smallskip

\begin{iii}
\item[iv) ]
In~this case, the assumption that the surface would have
Picard~rank~$2$
over~$\overline{\bbQ}$
implies that the discriminants of both Picard groups,
$\Pic(S_{\overline{\bbQ}})$
and~$\Pic(S_{\overline{\bbF}_{\!p}})$,
are in the same square~class. Note~here that reduction
modulo~$p$
respects the intersection~product.
\item[v) ]
To~obtain a contradiction, one combines information from two~primes. It~may happen that one has a rank bound
of~$2$
at both places but different
square classes for the discriminant do~arise. Then,~these data are incompatible
with Picard
rank~$2$
over~$\overline{\bbQ}$.
\end{iii}
\end{meto}

\begin{imp}
Approach~\ref{vLm} accepts the possibility that
$\smash{\Pic(S_{\overline{\bbQ}}) \subset \Pic(S_{\overline{\bbF}_{\!p}})}$
might be a proper sublattice of full~rank. If~that occurred then one knows at~least that the two discriminants differ by a perfect~square. This~is a standard observation from the theory of~lattices.

We~will show in this article that such provisions need not be~made. From~the technical point of view, our main result states that, at least
for~$p \neq 2$,
the quotient
$\Pic(S_{\overline{\bbF}_{\!p}})/\Pic(S_{\overline{\bbQ}})$
is always torsion-free. This~is true actually in much more generality than just for
$K3$~surfaces.
It~follows in a rather straightforward manner from deformation theory, a tool developed by A.~Grothendieck and M.~Artin in the sixties of the last~century. To~be precise, our result is as~follows.
\end{imp}

\begin{theo}
\label{maini}
Let\/~$p\neq 2$
be a prime number and\/
$X$
be a scheme proper and flat
over\/~$\bbZ$.
Suppose~that the special
fiber\/~$X_p$
is~non-singular and satisfies
$H^1(X_p, \calO_{X_p}) = 0$.\smallskip

\noindent
Then,~the specialization homomorphism\/
$\Pic(X_{\overline\bbQ}) \to \Pic(X_{\overline{\bbF}_{\!p}})$
has a torsion-free~cokernel.
\end{theo}

\begin{rems}
\begin{abc}
\item
Recall~that, for a
$K3$~surface~$S$,
one
has~$H^1(S,\calO_S) = 0$
\cite[Chap.~VI, Table~10]{BPV}.
\item
We~will prove this theorem in~\ref{main}. As~an application, one may prove
$\rk \Pic(S_{\overline\bbQ}) = 1$
for
a~$K3$~surface~$S$
using its reduction only at a single~prime.
This~works as~follows.
\end{abc}\bigskip

\noindent
{\bf Approach.}
Let
a~$K3$~surface~$S$
be~given.

\begin{iii}
\item
For~a
prime~$p \neq 2$
of good reduction, perform steps~i), ii) and~iii) as in~\ref{vLm}. Thereby,~the hope is to prove
$\rk \Pic(S_{\overline\bbF_{\!p}}) \leq 2$.
Further,~compute the discriminant giving two explicit~generators. Alternatively,~one might use the Artin-Tate~formula.
\item
Assume~$\rk \Pic(S_{\overline\bbQ}) = 2$.
Then,~according to Theorem~\ref{maini}, every invertible sheaf
on~$S_{\overline\bbF_{\!p}}$
lifts
to~$S_{\overline\bbQ}$.
Use~reduction theory of binary quadratic forms or explicit arguments to estimate the degree of a hypothetical effective~divisor. Finally,~use Gr\"obner bases to verify that such a divisor does not~exist.
\end{iii}
\end{rems}

\begin{ex}
Consider the
$K3$~surface~$S$
over~$\bbQ$,
given~by
$$w^2 = x^5 y + x^4 y^2 + 2 x^3 y^3 + x^2 y^4 + x y^5 + 4 y^6 + 2 x^5 z + 2 x^4 z^2 + 4 x^3 z^3 + 2 x z^5 + 4 z^6 \, .$$
Then,~$\rk \Pic(S_{\overline\bbQ}) = 1$.\smallskip

\noindent
{\bf Proof.}
For~the reduction
of~$S$
at the
prime~$5$,
one sees that the branch locus has a tritangent line given
by~$z - 2y = 0$.
It~meets the branch locus at
$(1:0:0)$,
$(1:3:1)$,
and~$(0:1:2)$.

The~numbers of points
over~$\bbF_{5^d}$
are, in this order,
$41$,
$751$,
$15\,626$,
$392\,251$,
$9\,759\,376$,
$244\,134\,376$,
$6\,103\,312\,501$,
$152\,589\,156\,251$,
$3\,814\,704\,296\,876$,
and
$95\,367\,474\,609\,376$.
Thus,~the traces of the Frobenius
on~$\smash{H^2_\et (S_{\overline\bbF_5} \!, \overline\bbQ_l)}$
are
$15$,
$125$,
$0$,
$1\,625$,
$-6\,250$,
\,$-6\,250$,
\,$-203\,125$,
\,$1\,265\,625$,
\,$7\,031\,250$,
and~$42\,968\,750$.
Algorithm~23 of~\cite{EJ} shows that the sign in the functional equation is~positive. The~characteristic polynomial of the Frobenius is therefore completely~determined. For~its decomposition into prime polynomials, we find (after~scaling)
\begin{eqnarray*}
\begin{aligned}
& (t-5)^2
(t^{20} - 5\,t^{19} - 25\,t^{18} + 250\,t^{17} - 250\,t^{16} - 1\,875\,t^{15} + 12\,500\,t^{14} - 31\,250\,t^{13} \\
& {} \hspace{1.8cm} - 156\,250\,t^{12} + 390\,625\,t^{11} + 5\,859\,375\,t^{10} + 9\,765\,625\,t^9 - 97\,656\,250\,t^8 \\
& {} \hspace{1.5cm} - 488\,281\,250\,t^7 + 4\,882\,812\,500\,t^6 - 18\,310\,546\,875\,t^5 - 61\,035\,156\,250\,t^4 \\
& {} \hspace{2.6cm} + 1\,525\,878\,906\,250\,t^3 - 3\,814\,697\,265\,625\,t^2 - 19\,073\,486\,328\,125\,t \\
& {} \hspace{10.5cm} + 95\,367\,431\,640\,625) \, .
\end{aligned}
\end{eqnarray*}
This~shows~$\rk \Pic(S_{\overline\bbF_{\!5}}) \leq 2$.

The~splits of the pull-back of the tritangent line are explicit generators
for~$\Pic(S_{\overline\bbF_{\!5}})$.
Such~a
split~$l$,
being a projective line, has self-intersection number
$l^2 = -2$.
Further,~$lh = 1$
for
$h$
the pull-back of a~line. If~we had
$\rk \Pic(S_{\overline\bbQ}) = 2$
then the invertible
sheaf~$\calO(l)$
would lift
to~$S_{\overline\bbQ}$.
We~had a
divisor~$L$
on~$S_{\overline\bbQ}$
such that
$HL = 1$
and~$L^2 = -2$.
By~\cite[Ch.\,VIII, Proposition~3.6.i]{BPV}, such a divisor is automatically~effective.

$HL = 1$~shows
that~$L$
is obtained from a line
on~$\Pb^2$,
the pull-back of which splits into two~components. This~is possible only for a tritangent line of the branch~locus. \cite[Algorithm~8]{EJ}~shows, however, that such a tritangent line does not~exist.
\eop
\end{ex}

\section{The sequence of the Picard lattices}

\begin{remark}
The~proof of Theorem~\ref{maini} relies on deformation-the\-o\-retic methods~\cite{Ar,Kl}. For~$K3$~surfaces
and
prime-to-$p$~torsion,
one could have used \'etale cohomology which appears to be more~natural.

In~fact, to show
$\smash{\Pic(X_{\overline{\bbF}_{\!p}})/\Pic(X_{\overline\bbQ})}$
has no
\mbox{$l$-torsion},
it is sufficient to consider
$\smash{\Pic(X_{\overline{\bbF}_{\!p}}) \!\otimes_\bbZ\! \bbZ_l/\Pic(X_{\overline\bbQ}) \!\otimes_\bbZ\! \bbZ_l}$.
But~$\smash{\Pic(X_{\overline{\bbF}_{\!p}}) \!\otimes_\bbZ\! \bbZ_l \subseteq H^2_\et(X_{\overline{\bbF}_{\!p}}, \bbZ_l(1))}$
which, by standard comparison theorems, is isomorphic
to~$H^2_\sing(X(\bbC), \bbZ) \!\otimes_\bbZ\! \bbZ_l$.
On~the other hand,
$\Pic(X_{\overline\bbQ}) \cong \Pic(X_\bbC)$.
Finally,~$H^2_\sing(X(\bbC), \bbZ) / \Pic(X_\bbC)$
is torsion-free according to the Lefschetz
$(1,1)$-theorem.
\end{remark}

\begin{nota}
Let~$X$
be a
\mbox{$\bbZ_p$-scheme}.
Then,~we will write
$X_p$
for the special fiber and, more generally,
$X_{p^n} := X \!\times_{\Spec \bbZ_p}\! \Spec \bbZ/p^n\bbZ$.
Finally,~let
$\smash{\widehat{X}}$
be the formal~scheme obtained by
completing~$X$
along~$(p)$.
\end{nota}

\begin{lem}
\label{eins}
Let\/~$p\neq 2$
be a prime number and\/
$X$
a\/
$\bbZ_p$-scheme
which is Noetherian,separated, and
fulfills\/~$H^1(X_p, \calO_{X_p}) = 0$.
Denote~by\/
$P \subseteq \Pic(X_p)$
the subset of all invertible sheaves allowing a lift as an invertible sheaf
on\/~$\smash{\widehat{X}}$.\smallskip

\noindent
Then,
$\Pic(X_p)/P$
is torsion-free.\medskip

\noindent
{\bf Proof.}
{\em
{\em First step.}
Preliminaries.\smallskip

\noindent
Assume,~to the contrary, that
$\Pic(X_p)/P$
has~torsion. Then,~there are a prime
number~$l$
and an invertible
sheaf~$\calL \in \Pic(X_p) \!\setminus\! P$
such
that~$\calL^{\otimes l} \in P$.
This~means that
$\calL^{\otimes l}$
lifts
to~$\smash{\widehat{X}}$
but
$\calL$
does~not. We~have to show that this situation is~impossible.

By~\cite[Proposition~II.9.6]{Ha}, an invertible sheaf
on~$\smash{\widehat{X}}$
is the same as an inverse system
$(\calI_n)_n$
of invertible sheaves
$\calI_n \in \Pic(X_{p^n})$
such that
$\calI_{n+1} |_{X_{p^n}} = \calI_n$
for~all~$n$.
By~assumption, we have such a system
for~$\calI_0 = \calL^{\otimes l}$.
It~has to be shown that the invertible sheaf
$\calL$,
too, lifts
to~$X_{p^n}$
for
all~$n$.\medskip

\noindent
{\em Second step.}
Obstructions.\smallskip

\noindent
We~will construct
sheaves~$\calL_n \in \Pic(X_{p^n})$
lifting~$\calL$,
inductively. These~will satisfy, in addition, the
relation~$\calL_n^{\otimes l} \cong \calI_n$.
First,~we
put~$\calL_0 := \calL$.

For~the induction step, consider the short exact~sequence
$$0 \longrightarrow \calK \longrightarrow \calO^*_{X_{p^{n+1}}} \longrightarrow \calO^*_{X_{p^n}} \longrightarrow 0 \, .$$
Here,~we have
$\calO_{X_p} \cong \calK$
via the exponential map
$x \mapsto 1 + p^n x \pmod {p^{n+1}}$.
This~yields the commutative diagram with exact~rows,
$$
\diagram
0 \rto & \Pic(X_{p^{n+1}}) \rto\dto^{(.)^{\otimes l}} & \Pic(X_{p^n}) \rto\dto^{(.)^{\otimes l}} & H^2(X_p, \calO_{X_p}) \phantom{\, .} \dto^{\cdot l} \\
0 \rto & \Pic(X_{p^{n+1}}) \rto & \Pic(X_{p^n}) \rto & H^2(X_p, \calO_{X_p}) \, . 
\enddiagram
$$
The~group~$H^2(X_p, \calO_{X_p})$
is
\mbox{$p$-torsion}
as the sheaf
$\calO_{X_p}$
is annihilated
by~$p$.
In~particular, it is uniquely
\mbox{$l$-divisible}.
Further,
$\calI_n \in \Pic(X_{p^n})$
is the image of
$\calI_{n+1} \in \Pic(X_{p^{n+1}})$
and~$\calL_n \in \Pic(X_{p^n})$.
A~standard diagram argument
yields some invertible
sheaf~$\calL_{n+1} \in \Pic(X_{p^{n+1}})$
which is mapped to
$\calL_n$
and~$\calI_{n+1}$.
This~completes the proof
for~$l \neq p$.\medskip

\noindent
{\em Third step.}
The case~$l=p$.\smallskip

\noindent
Here,~we first observe the congruence
$$(1 + p^n c)^p \equiv 1 + p^{n+1} c \pmod {p^{n+2}}$$
which,
as~$p > 2$,
is valid for
every~$n \geq 1$.
This~has the striking consequence that,
for~$s \in \Gamma(U, \calO_{X_{p^n}}^*)$,
the
power~$s^p$
automatically defines a section
of~$\calO_{X_{p^{n+1}}}^*$.
Further,~we have the commutative diagram
$$
\definemorphism{gleich}\Solid\notip\notip
\diagram
0 \rto& \calO_{X_p} \rto\dgleich& \calO^*_{X_{p^{n+1}}} \rto\dto^{(.)^p}& \calO^*_{X_{p^n}} \rto\dto^{(.)^p}& 0 \\
0 \rto& \calO_{X_p} \rto& \calO^*_{X_{p^{n+2}}} \rto& \calO^*_{X_{p^{n+1}}} \rto& 0
\enddiagram
$$
with exact~rows. Taking~cohomology, this yields the commutative diagram with exact~rows,
$$
\definemorphism{gleich}\Solid\notip\notip
\diagram
0 \rto & \Pic(X_{p^{n+1}}) \rto\dto^{(.)^{\otimes p}} & \Pic(X_{p^n}) \rto\dto^{(.)^{\otimes p}} & H^2(X_p, \calO_{X_p}) \phantom{\, .} \dgleich \\
0 \rto & \Pic(X_{p^{n+2}}) \rto & \Pic(X_{p^{n+1}}) \rto & H^2(X_p, \calO_{X_p}) \, . 
\enddiagram
$$
We~see, in particular, that the lift of an invertible sheaf, if possible, is unique up to~isomorphism.

We~will inductively construct a sequence of
sheaves~$\calL_n \in \Pic(X_{p^n})$
lifting~$\calL$
such
that~$\calL_n^{\otimes p} \cong \calI_n$.
To~start, simply
put~$\calL_0 := \calL$.
For~the induction step, we observe that
$\calL_n^{\otimes p} \cong \calI_n$
implies that
$\calL_n$
is mapped to
$\calI_{n+1}$
under the middle vertical arrow in the~diagram. Indeed,~the lifting of an invertible sheaf is~unique. The~same diagram argument as in the second step completes the~proof.
}
\eop
\end{lem}

\begin{remark}
For~$p=2$,
the same argument shows that
$\Pic(X_2)/P$
may only have
\mbox{$2$-power}~torsion.
\end {remark}

\begin{ttt}
To~illustrate the effect of the obstructions, suppose that
$\Pic(X_p) = \bbZ^n$
and~$H^2(X_p, \calO_{X_p}) \cong \bbF_{\!p}$.
Then,~the lattices
$\Lambda_i := \Pic(X_{p^i})$
form a
system~$\{\Lambda_i\}_{i \in \bbN}$
such~that
$$\ldots \subseteq \Lambda_i \subseteq \ldots \subseteq \Lambda_2 \subseteq \Lambda_1 \, ,$$
$\Lambda_i \subseteq \Lambda_{i-1}$
is always of
index~$1$
or~$p$,
and~$\calL^{\otimes p} \in \Lambda_i$
if and only
if~$\calL \in \Lambda_{i-1}$.

According~to Lemma~\ref{Stephan} below, the system
$\{\Lambda_i \!\otimes_\bbZ\! \bbZ_p\}_{i \in \bbN}$
is isomorphic to
$\{\bbZ_p \!\oplus \cdots \oplus\! \bbZ_p \!\oplus\! p^i \bbZ_p\}_{i \in \bbN}$.
I.e.,~there is a linear functional
$$H\colon \Pic(X_p) = \bbZ^n \to \bbZ_p, \quad (x_1,\ldots,x_n) \mapsto a_1x_1 + \cdots +a_nx_n$$
with coefficients
$a_1,\ldots,a_n \in \bbZ_p$
such that, for
$\calL \in \Pic(X_p)$
arbitrary,
$p^i | H(\calL)$
if and only if
$\calL$
lifts
to~$\Pic(X_{p^i})$.

$H$~somehow
collects all the obstruction maps into a single~homomorphism.
Further,~$H(\calL) = 0$
if and only
if~$\calL \in P$.
This~shows again
that~$\Pic(X_p)/P \hookrightarrow \bbZ_p$
is torsion-free.
\end{ttt}

\begin{remark}
This~formulation also indicates that it is difficult to show
$\rk P \leq \rk \Pic(X_p) - 2$.
For~this, one had to ensure that the
\mbox{$\bbZ$-rank}
of~$\im H$
is at
least~$2$.
But~this is impossible knowing only
\mbox{$p$-adic}~approximations
of~$a_1, \ldots, a_n$.
\end{remark}

\begin{lem}
\label{Stephan}
Let\/~$\{\Lambda_i\}_{i\in\bbN}$
be a sequence of\/
\mbox{$p$-adic}
lattices such~that

\begin{iii}
\item
$\Lambda_{i+1} \subset \Lambda_i$,
\item
$\Lambda_i / \Lambda_{i+1} = \bbZ/p\bbZ$,
\item
$x \in \Lambda_i \!\setminus\! \Lambda_{i+1} \Longrightarrow px \in \Lambda_{i+1} \!\setminus\! \Lambda_{i+2}$.
\end{iii}
Then,~there exists a basis\/
$(b_1, \ldots, b_n)$
of\/~$\Lambda_1$
such that\/
$\Lambda_i = \langle b_1, \ldots, b_{n-1}, p^{i-1} b_n \rangle$.\smallskip

\noindent
{\bf Proof} {\em (cf.\,\cite{We})}{\bf .}
{\em
We~first observe
that~$\Lambda_1/\Lambda_i \cong \bbZ/p^{i-1}\bbZ$.
Indeed,~the quotient
$\Lambda_1 / \Lambda_i$
is precisely of
order~$p^{i-1}$.
Further,~for
$x \in \Lambda_1 \!\setminus\! \Lambda_2$,
we find
$px \in \Lambda_2 \!\setminus\! \Lambda_3$
and, finally,
$p^{i-2} x \in \Lambda_{i-1} \!\setminus\! \Lambda_i$.
In~particular, we see that
$\Lambda_0 / \Lambda_i$
has an element of
order~$p^{i-1}$.

Let~now
$i$
be~fixed. By~the elementary divisor theorem, there exists a basis
$(b_1, \ldots, b_n)$
of~$\Lambda_0$
such that
$(p^{e_1} b_1, \ldots, p^{e_n} b_n)$
is a basis
of~$\Lambda_i$.
As~this yields
$\Lambda_1/\Lambda_i \cong \bbZ /p^{e_1}\bbZ \times \cdots \times \bbZ/p^{e_n}\bbZ$,
we may conclude
$e_1 = \cdots = e_{n-1} = 0$
and~$e_n = i-1$.
The~only lattices between 
$\Lambda_0$
and~$\Lambda_i$
are
$\langle b_1, \ldots, b_{n-1}, p^j b_n \rangle$
for~$j = 1, \ldots, i-2$.
Thus,~we have shown the assertion for a finite chain of~lattices.

To~prove it for the infinite sequence, we observe that the space of all bases
of~$\Lambda_1$
is compact in the
\mbox{$p$-adic}~topology.
For~every
$i \in \bbN$,
there is a basis
$B^{(i)} = (b^{(i)}_1, \ldots, b^{(i)}_n)$
of~$\Lambda_1$
having the desired property for the finite subsequence
$\Lambda_1, \ldots ,\Lambda_i$.
Consider~the limit
$(b_1, \ldots ,b_n)$
of a convergent subsequence
of~$\{B^{(i)}\}_{i\in\bbN}$.

We~claim that
$(b_1, \ldots, b_{n-1}, p^{i-1} b_n)$
is a basis
for~$\Lambda_i$.
Indeed,~$(b_1, \ldots, b_{n-1}, p^{i-1} b_n)$
is arbitrarily close to a basis which completes the~proof.
}
\eop
\end{lem}

\section{The quotient
\boldmath$\Picb(X_{\overline{\bbF}_{\!p}})/\Picb(X_{\overline\bbQ})$}

\begin{subl}
\label{Zariski}
Let\/~$p$
be a prime number and\/
$X$
be a\/
\mbox{$\bbZ_p$-scheme}
which is proper and~flat. Suppose~that the generic
fiber\/~$X_\eta$
is connected and the special
fiber\/~$X_p$
is~non-singular.\smallskip

\noindent
Then,~$X_p$
is~irreducible.\smallskip

\noindent
{\bf Proof.}
{\em
The~function field
$K := \Gamma(X_\eta, \calO_{X_\eta})$
is a finite extension
of~$\bbQ_p$.
Further,~$O := \Gamma(X, \calO_X)$~is
a finite
\mbox{$\bbZ_p$-algebra}
being an integral domain with quotient
field~$K$.
Clearly,~$O/pO$
is contained
in~$\Gamma(X_p, \calO_{X_p})$.
But,~according to the assumption, the latter does not have nilpotent elements other than~zero.
Hence,~$p$
generates the maximal ideal
of~$O$.
This~means,
$K/\bbQ_p$
is necessarily unramified
and~$O = \calO_K$
is its ring of~integers. Stein~factorization provides us with a morphism
$X \to \Spec \calO_K$
with connected~fibers. From~this, we immediately see that
$X_p$
is~connected.
As~$X_p$
is non-singular, this is enough for~irreducibility.
}
\eop
\end{subl}\pagebreak[3]

\begin{lem}
\label{zwei}
Let\/~$p\neq 2$
be a prime number and\/
$X$
be a\/
\mbox{$\bbZ_p$-scheme}
which is proper and~flat. Suppose~that the special
fiber\/~$X_p$
is~non-singular and satisfies
$H^1(X_p, \calO_{X_p}) = 0$.\smallskip

\noindent
Then,~the specialization homomorphism\/
$\sp\colon \Pic(X_\eta) \to \Pic(X_p)$
from the generic fiber has a torsion-free~cokernel.\smallskip

\noindent
{\bf Proof.}
{\em
As~each connected component may be treated separately, we assume without restriction that
$X$~is~connected.
Further,~the assumption implies that
$X$
is non-singular.
Hence,~$X$
is actually~irreducible. This~implies that
$X_\eta$
is irreducible,~too. Finally,~we conclude irreducibility
of~$X_p$
from Sublemma~\ref{Zariski}.

There~is a specialization map
$\Pic(X_\eta) \to \Pic(X)$
given by taking the Zariski closure
in~$X$
of a Weil divisor
on~$X_\eta$.
This~map is injective as the restriction forms a section to~it. It~is a surjection, too, as the only vertical divisors are principal, associated to the powers
of~$(p)$.

Further,~by A.~Grothendieck's existence theorem~\cite[Corollaire~(5.1.6)]{EGAIII}, one
has~$\smash{\Pic(X) = \Pic(\widehat{X})}$.
The~assertion now follows from Lemma~\ref{eins}.
}
\eop
\end{lem}

\begin{coro}
Let\/~$p\neq 2$
be a prime number and\/
$X$
be a\/
\mbox{$\bbZ_p$-scheme}
which is proper and~flat. Suppose~that the special
fiber\/~$X_p$
is~non-singular and satisfies
$H^1(X_p, \calO_{X_p}) = 0$.
Further,~let\/~$K/\bbQ_p$
be an unramified field extension and denote the residue field
of\/~$K$
by\/~$k$.\smallskip

\noindent
Then,~the cokernels of the specialization homomorphisms

\begin{iii}
\item
$\sp_K \colon \Pic(X_K) \to \Pic(X_k)$,
\item
$\sp_{\bbQ_p^\nr} \colon \Pic(X_{\bbQ_p^\nr}) \to \Pic(X_{\overline{\bbF}_{\!p}})$,
and
\item
$\sp_{\overline{\bbQ}_p} \colon \Pic(X_{\overline{\bbQ}_p}) \to \Pic(X_{\overline{\bbF}_{\!p}})$
\end{iii}
are torsion-free.\smallskip

\noindent
{\bf Proof.}
{\em
i)
Apply Lemma~\ref{zwei} to the fiber product
$X \!\times_{\Spec \bbZ_p}\! \Spec \calO_K$.\smallskip

\noindent
ii)
As~the filtered direct limit functor is exact, the desired cokernel is the same as
$$\limdir \coker (\sp_K \colon \Pic(X_K) \to \Pic(X_k))$$
where
$K$
is running over the unramified extensions
of~$\bbQ_p$
and~$k$
denotes the residue field
of~$K$.
As~all the cokernels are torsion-free, the assertion~follows.\smallskip

\noindent
iii)
We~claim that
$\smash{\sp_{\overline{\bbQ}_p}}$
has the same image
in~$\smash{\Pic(X_{\overline{\bbF}_{\!p}})}$
as~$\smash{\sp_{\bbQ_p^\nr}}$.
Let~$\smash{\calL \in \Pic(X_{\overline{\bbQ}_p})}$.
The~Galois group
$\smash{\Gamma := \Gal(\overline{\bbQ}_p/\bbQ_p^\nr)}$
sends
$\calL$
to a finite
orbit~$\smash{\{\calL_1, \ldots, \calL_m\}}$.
The~specializations
of~$\calL_1, \ldots, \calL_m$
in~$\smash{\Pic(X_{\overline{\bbF}_{\!p}})}$
are all the~same. Therefore,
$$m \!\cdot\! \sp_{\overline{\bbQ}_p} (\calL) = \sp_{\overline{\bbQ}_p} (\calL^{\otimes m}) = \sp_{\overline{\bbQ}_p} (\calL_1 \!\otimes \cdots \otimes\! \calL_m) = \sp_{\bbQ_p^\nr} (\calL_1 \!\otimes \cdots \otimes\! \calL_m)$$
since
$\calL_1 \!\otimes \cdots \otimes\! \calL_m$
is
\mbox{$\Gamma$-invariant}.
Hence,~$\smash{m \!\cdot\! \sp_{\overline{\bbQ}_p} (\calL) \in \im \sp_{\bbQ_p^\nr}}$.
As~$\smash{\sp_{\bbQ_p^\nr}}$
has a torsion-free cokernel, we see that
$\smash{\sp_{\overline{\bbQ}_p} (\calL) \in \im \sp_{\bbQ_p^\nr}}$,~too.
}
\eop
\end{coro}

\begin{theo}
\label{main}
Let\/~$p\neq 2$
be a prime number and\/
$X$
be a scheme proper and flat
over\/~$\bbZ$.
Suppose~that the special
fiber\/~$X_p$
is~non-singular and satisfies
$H^1(X_p, \calO_{X_p}) = 0$.\smallskip

\noindent
Then,~the specialization homomorphism\/
$\sp_{\overline\bbQ} \colon \Pic(X_{\overline\bbQ}) \to \Pic(X_{\overline{\bbF}_{\!p}})$
has a torsion-free~cokernel.\smallskip

\noindent
{\bf Proof.}
{\em
There~is a canonical injection
$\smash{\Pic(X_{\overline\bbQ}) \hookrightarrow \Pic(X_{\overline\bbQ_p})}$.
We~have to show that both Picard groups have the same image under specialization
to~$\smash{\Pic(X_{\overline{\bbF}_{\!p}})}$.

For~this, we switch at first to the
scheme~$X_Z$
for~$Z = \bbZ[\frac1m]$
where~$m$
is an integer divisible by all primes of bad reduction but not
by~$p$.
Again,~we may assume without restriction that
$X_Z$
is~connected. Further,~by construction,
$X_Z$
is non-singular and, therefore,~irreducible. By~virtue of Sublemma~\ref{Zariski}, all the special fibers
of~$X_Z$
are~irreducible.

According~to a theorem of Grothendieck (cf.~\cite[Theorem~4.8]{Kl}), the Picard scheme~$\Picb_{X_Z/Z}$
exists in this situation as a scheme, locally of finite type
over~$Z$.
This~means, we are given a morphism
$\smash{i\colon \Spec \overline{\bbQ}_p \to \Picb_{X_Z/Z}}$
and have to show that there is a morphism
$\Spec \overline{\bbQ} \to \Picb_{X_Z/Z}$
such that the specializations
modulo~$p$
are the~same.

Locally,~near the image
of~$i$,
we have an affine open subset
$U \cong \Spec R \subseteq \Picb_{X_Z/Z}$
for
$R$
a finitely generated
\mbox{$Z$-algebra}.
We~are thus given a ring homomorphism
$\iota \colon R \to \overline{\bbQ}_p$.
This~actually
maps~$R$
to
$\calO_K$
for a suitable finite
extension~$K/\bbQ_p$.
Unfortunately,~as
a~\mbox{$Z$-algebra},
$\calO_K$~is
not finitely~generated.
On~the other hand,
$\im \iota =: S \subset \calO_K$~is
clearly a finitely generated
\mbox{$Z$-algebra}.
We~fix a set of generators
$\{T_1,\ldots,T_n\}$
of~$S$.

Preserving~the induced homomorphism to
$\bbF_{\!q} := \calO_K/\frakm_K$,
our goal is to
replace~$\iota$
by a homomorphism to another subring
$S^\prime \subset \calO_K$
such that
$S$
is finite as
a~\mbox{$Z$-module}.
For~this, we will construct an algebra~homomorphism
$\varphi\colon S \to S^\prime$
such that
$\nu(x - \varphi(x)) \geq 1$
for
every~$x \in S$.
Here,~$\nu$
denotes the discrete valuation
on~$\calO_K$.

To~perform this construction, we apply Noether normalization~\cite[Ch.\,V, \S4, Theorem\,8]{ZS}
to~$S \otimes_\bbZ \bbQ$.
This~states that
$S \otimes_\bbZ \bbQ$~is
an~integral extension of a polynomial
ring~$\bbQ[X_1,\ldots,X_k] \subseteq S \!\otimes_\bbZ\! \bbQ$.
We~send
$X_1,\ldots,X_k$
to elements
of~$\calO_K$
algebraic
over~$\bbQ$
such
that~$\nu(X_i - \varphi(X_i)) \gg 0$.
Then,~this extends to a homomorphism of the whole
of~$S \!\otimes_\bbZ\! \bbQ$.
We~claim
that~$\smash{\nu(T_i - \varphi(T_i)) \geq 1}$
for~$i = 1, \ldots, n$.
Indeed,~as
$T_i$
is integral over
$\bbQ[X_1,\ldots,X_k,T_1,\ldots,T_{i-1}]$,
this follows from an iterated application of Hensel's lemma in the form
of~\cite[Proposition~5.5]{Na}.

Since~$S$
is generated
by~$T_1, \ldots, T_n$
as
a~\mbox{$Z$-algebra}
and
$\nu(z) \geq 0$
for every
$z \in Z$,
we see that
$\nu(x - \varphi(x)) \geq 1$
for
every~$x \in S$.
This~completes the~proof.
}
\eop
\end{theo}

\section{An explicit obstruction}

\begin{prop}
Let\/~$S$
be a\/
$K3$~surface
of
degree\/~$2$
over\/~$\bbQ$,
given explicitly~by
$$w^2 = f_6(x,y,z)$$
for\/~$f_6 \in \bbZ[x,y,z]$
of
degree\/~$6$.
Suppose,~for a prime\/
$p \neq 2$,
there is an\/
\mbox{$\bbF_{\!p}$-rational}
tritangent line\/
``$\ell = 0$\!''
of the ramification locus
of\/~$S_p$.
Write\/~$l$
for a split of the pull-back of the~tritangent.\smallskip

\noindent
One~has\/
$f_6 \equiv f_3^2 + \ell f_5 \pmod p$
for homogeneous forms\/
$f_3, f_5 \in \bbZ[x,y,z]$.
Put
$$G(x,y,z) := (f_6 - f_3^2 - \ell f_5) / p \, .$$
Then,~the obstruction to lifting\/
$\calO(l)$
to\/~$S_{p^2}$
is\/~$((-G) \!\mod (p,\ell,f_3,f_5))$.\medskip

\noindent
{\bf Proof.}
{\em
{\em First step.}
An affine open covering
of~$S_p$.\smallskip

\noindent
On~$S_p$,
we have
$w^2 = f_3^2 + \ell f_5$
and,
for~$h$
a quadric,
$w^2 = (f_3 + \ell h)^2 + \ell f_5^\prime$
where
$f_5^\prime := f_5 - 2f_3h - \ell h^2$.
On~``$\ell = 0$'',
$f_3$
and~$f_5$
have no common zero as this would cause a singularity
on~$S_p$.
Hence,~for a suitably
chosen~$h$,
the three forms
$\ell$,
$f_5$,
and~$f_5^\prime$
do not have a common~zero. For~this, it may be necessary to extend the ground~field.
The~sets
``$\ell \neq 0$'',
``$f_5 \neq 0$'',
and~``$f_5^\prime \neq 0$'',
form an affine open covering
of~$S_p$.
We~may extend them in the obvious manner to an affine open covering
of~$S$.\medskip

\noindent
{\em Second step.}
The invertible
sheaf~$\calO(5l)$.\smallskip

\noindent
We~start with
$\calO(5l)$
instead
of~$\calO(l)$
as this will turn out to be~easier.
$\calO(5l)$~is
given by the rational functions
$1$
on~``$\ell \neq 0$'',
$\smash{\frac{f_5^3}{(w+f_3)^5}}$
on~``$f_5 \neq 0$'',
and
$\smash{\frac{f_5^\prime {}^3}{(w+f_3+\ell h)^5}}$
on~``$f_5^\prime \neq 0$''.
Thus,~the transition functions are
\begin{eqnarray*}
& & \frac{f_5^3}{(w+f_3)^5} = \frac{(w-f_3)^5}{\ell^5 f_5^2} \, ,
\qquad \frac{(w+f_3)^5 f_5^\prime {}^3}{(w+f_3+\ell h)^5 f_5^3} = \frac{(w+f_3)^5 (w-f_3-\ell h)^5}{\ell^5 f_5^3 f_5^\prime {}^2} \, , \\
& & {\rm and~} \quad \frac{(w+f_3+\ell h)^5}{f_5^\prime {}^3} \, .
\end{eqnarray*}\pagebreak[3]

\noindent
{\em Third step.}
The~obstruction.\smallskip

\noindent
We~may lift the first and third transition functions~naively. The~middle one is a transition function between
``$f_5 \neq 0$''
and~``$f_5^\prime \neq 0$''
and, thus, must not have a pole
at~``$\ell \neq 0$''.
We~lift
$(w+f_3) (w-f_3)$
as~$\ell f_5$
and obtain, in total,
$$\frac{[f_5 - h(w + f_3)]^5}{f_5^3 f_5^\prime {}^2} \, .$$
The~product of the three lifts~is
$$\frac{(w-f_3)^5 [f_5 - h(w + f_3)]^5 (w+f_3+\ell h)^5}{\ell^5 f_5^5 f_5^\prime {}^5} \, .$$
Observe~that, in the form described, the transition functions may be lifted even to the affine~open subsets
of~$S$,
not just
to~$\smash{S_{p^2}}$.
Hence,~the exponential of the obstruction
for~$\calO(l)$
is~$\smash{\frac{(w-f_3) [f_5 - h(w + f_3)] (w+f_3+\ell h)}{\ell f_5 f_5^\prime}}$,
also in the case
that~$p = 5$.

Evaluating~this expression, making use of the identity
$w^2 - f_3^2 = \ell f_5 + pG$,
we end up~with
$\smash{1 + p \,\frac{G(f_5 - hw - hf_3 - \ell h^2)}{\ell f_5 f_5^\prime}}$.
Therefore,~the obstruction to
lifting~$\calO(l)$
is given by the
\v{C}ech~cocycle
$$\frac{G(f_5 - hw - hf_3 - \ell h^2)}{\ell f_5 f_5^\prime} \, .$$\vskip\medskipamount

\noindent
{\em Fourth step.}
Simplification.\smallskip

\noindent
Any~rational function having poles in only two of the three divisors considered is a \v{C}ech~coboundary. Without~changing the cohomology class, we may therefore add to the numerator forms being homogeneous of
degree~$11$
and belonging to the
ideal~$(\ell,f_5,f_5^\prime)$.

On~the
line~``$\ell=0$'',
$f_5$
and~$f_5^\prime$
have no zeroes in~common. Thus,~they are coprime in the graded
ring~$\bbF_{\!p}[x,y,z]/(\ell)$.
Consequently,~$f_5$
and~$f_5^\prime$
already generate the full
\mbox{$10$-dimensional}
space of forms of
degree~$9$.
Even~more, they must generate the space of forms of
degree~$11$.
This~shows that we may simplify the \v{C}ech cocycle
to~$\smash{\frac{-Ghw}{\ell f_5 f_5^\prime}}$.

Hence,~the obstruction to
lifting~$\calO(l)$
is~$((-Gh) \!\mod (p,\ell,f_5,f_5^\prime))$.
The~ideal is the same
as~$(p,\ell,hf_3,f_5)$.
Thus,~the question is whether
\mbox{$((-Gh) \!\mod (p,\ell))$}
is a combination
of~$hf_3$
and~$f_5$.
As,~on the
line~``$\ell=0$''
on~$S_p$,
$h$
and~$f_5$
have no common zeroes, they are~coprime.
$((-Gh) \!\mod (p,\ell))$
must be a combination
of~$hf_3$
and~$hf_5$.
We~may, as well,
consider~$((-G) \!\mod (p,\ell,f_3,f_5))$.
}
\eop
\end{prop}

\begin{ex}
Let~$S$
be a
$K3$~surface
over~$\bbQ$
given by
$w^2 = f_6(x,y,z)$.
Suppose
\begin{eqnarray*}
f_6(x,y,z) & \equiv &
x^6 + 2x^5z + 2x^4y^2 + 2x^4z^2 + 2x^3y^3 + 2x^3z^3 \\
 & & \hspace{2cm}{} + 2x^2y^4 + 2x^2y^3z + x^2z^4 + xy^3z^2 + 2xz^5 + y^6 \pmod 3 \, .
\end{eqnarray*}
Assume~further that the coefficient
of~$y^2 z^4$
is not divisible
by~$9$.\smallskip

\noindent
Then,~$\rk \Pic(S_{\overline\bbQ}) = 1$.\smallskip

\noindent
{\bf Proof.}
A~direct calculation shows that,
modulo~$3$,
the right hand side
is~$f_3^2 + xf_5$
for
$f_3 = 2x^3 + 2x^2z + xz^2 + 2y^3$
and
$f_5 = 2x^3y^2 + x^2z^3 + 2xy^4 + 2z^5$.
Thus,~the branch locus
of~$S_3$
has a tritangent line given
by~$x = 0$.

The~numbers of points
over~$\bbF_{\!3^d}$
are, in this order,
$19$,
$127$,
$676$,
$6\,751$,
$58\,564$,
$532\,414$,
$4\,791\,232$,
$43\,038\,703$,
$387\,383\,311$,
and~$3\,486\,675\,052$.
For~the decomposition of the characteristic polynomial of the Frobenius, we~find
\begin{eqnarray*}
\begin{aligned}
& (t-3)^2
(t^{20} - 3\,t^{19} - 9\,t^{18} + 72\,t^{17} - 8\,1t^{16} - 324\,t^{15} + 1\,458\,t^{14} - 2\,916\,t^{13} \\
& {} \hspace{0.7cm} + 4\,374\,t^{12} + 26\,244\,t^{11} - 137\,781\,t^{10} + 236\,196\,t^9 + 354\,294\,t^8 - 2\,125\,764\,t^7 \\
& {} \hspace{0.4cm} + 9\,565\,938\,t^6 - 19\,131\,876\,t^5 - 43\,046\,721\,t^4 + 344\,373\,768\,t^3 - 387\,420\,489\,t^2 \\
& {} \hspace{8.2cm} - 1\,162\,261\,467\,t + 348\,6784\,401) \, .
\end{aligned}
\end{eqnarray*}
This~shows~$\rk \Pic(S_{\overline\bbF_{\!3}}) \leq 2$.

Let~$l$
be a split of the pull-back of the tritangent~line. We~have to show that the obstruction to
lifting~$\calO(l)$
is non-zero. For~this, we observe that
$x$,
$f_3$,
and~$f_5$
do not generate the
monomial~$y^2z^4$.
However,~$G$
contains this monomial by its very~definition.
\eop
\end{ex}

\begin{ex}
Consider the
$K3$~surface~$S$
over~$\bbQ$,
given by
$w^2 = f_6(x,y,z)$
for
\begin{eqnarray*}
f_6(x,y,z) & = & 4x^6 + 2x^5y + 12x^5z + 2x^4y^2 + 4x^4yz + 12x^4z^2
+ 24x^3y^3 - 57x^3y^2z \\
 & & \hspace{0.15cm} {} - 9x^3yz^2 + 6x^3z^3 + 8x^2y^4 - 5x^2y^3z
- 72x^2y^2z^2 + 7x^2yz^3 + 4x^2z^4 \\
 & & \hspace{0.3cm} {} + 20xy^4z - 52xy^3z^2 - 57xy^2z^3
+ 7xyz^4 + 4y^5z - 7y^4z^2 - 18y^3z^3 \\
 & & \hspace{7.9cm} {} + 7y^2z^4 + 12yz^5 + 2z^6\, .
\end{eqnarray*}
Then,~$\rk \Pic(S_{\overline\bbQ}) = 3$.\smallskip

\noindent
{\bf Proof.}
We have
\begin{eqnarray*}
f_6 & = & (2x^3 + 2x^2z + 2y^2z + yz^2 + z^3)^2 \\
 & & {} + (2x^2 + 2xz + yz + z^2)(x^3y + 2x^3z + x^2y^2 + x^2yz + 2x^2z^2 + 12xy^3 \\
 & & {} \hspace{2.2cm} - 34xy^2z - 9xyz^2 - 2xz^3 + 4y^4 - 15y^3z - 7y^2z^2 + 9yz^3 + z^4)
\end{eqnarray*}
and
\begin{eqnarray*}
f_6 & = & 4(x^3 + 2x^2y + 2x^2z + xy^2 + xyz + xz^2 + y^2z + yz^2 + z^3)^2 \\
 & & {} - (x^2 + xz + yz + z^2)(14x^3y + 4x^3z + 22x^2y^2 + 22x^2yz + 8x^2z^2 - 8xy^3 \\
 & & {} \hspace{2.0cm} + 61xy^2z + 9xyz^2 + 6xz^3 - 4y^4 +  15y^3z + 11y^2z^2 - 6yz^3 + 2z^4) \, .
\end{eqnarray*}
Hence,~there are two conics
$C_1$
and~$C_2$
each of which is six times tangent to the ramification locus
of~$S$.
The~splits of their pull-backs yield the intersection~matrix
$$
\left(
\begin{array}{rrrr}
-2 &  6 &  1 &  3 \\
 6 & -2 &  3 &  1 \\
 1 &  3 & -2 &  6 \\
 3 &  1 &  6 & -2
\end{array}
\right)
$$
which is of
rank~$3$.
Hence,~$\rk \Pic(S_{\overline\bbQ}) \geq 3$.

On~the other hand,
$S$
has good reduction at the
prime~$p = 3$.
Point~counting over extensions
of~$\bbF_{\!3}$
shows that the characteristic polynomial of the Frobenius operating
on~$S_3$
is
\begin{eqnarray*}
 && {} (t-3)^4 (t^{18} + 3\,t^{17} + 6\,t^{16} + 18\,t^{15} + 108\,t^{14} + 405\,t^{13} + 972\,t^{12} + 2\,187\,t^{11} \\
 && {} \hspace{0.8cm} + 13\,122\,t^{10} + 52\,488\,t^9 + 118\,098\,t^8 + 177\,147\,t^7 + 708\,588\,t^6 + 2\,657\,205\,t^5 \\
 && {} \hspace{1.0cm} + 6\,377\,292\,t^4 + 9\,565\,938\,t^3 + 28\,697\,814\,t^2 + 129\,140\,163\,t + 387\,420\,489) \, .
\end{eqnarray*}
Consequently,~we
have~$\rk \Pic(S_{\overline\bbF_{\!3}}) \leq 4$.

In~particular, the assumption
$\rk \Pic(S_{\overline\bbQ}) > 3$
implies
$\rk \Pic(S_{\overline\bbQ}) = \rk \Pic(S_{\overline\bbF_{\!3}})$.
Theorem~\ref{main} guarantees that the specialization~map
$\smash{\sp_{\overline\bbQ} \colon \Pic(S_{\overline\bbQ}) \to \Pic(S_{\overline{\bbF}_{\!3}})}$
must be~bijective. Giving~one line bundle
$\calL \in \Pic(S_{\overline{\bbF}_{\!3}})$
with a non-trivial obstruction will be enough to yield a~contradiction.

For~this, observe that the ramification locus
of~$S_3$
has a tritangent line given
by~$x + y + z = 0$.
Indeed,
\begin{eqnarray*}
f_6(x,y,z) & \equiv & (x^3 + x^2y + xy^2 + y^3)^2 + (x+y+z)(2x^3y^2 + x^3yz + 2x^2yz^2 + 2xy^4 \\
 & & {} + xy^3z + xy^2z^2 + 2xyz^3 + xz^4 + 2y^5 + 2y^4z + yz^4 + 2z^5) \pmod 3 \, .
\end{eqnarray*}
Modulo~the
ideal~$(3,x+y+z)$,
we have
$f_3 \equiv x^3 + x^2y + xy^2 + y^3$,
$f_5 \equiv -(x^5 + x^3y^2 + x^2y^3 + xy^4 + y^5)$,
and~$G \equiv x^6 + 2x^5y + x^4y^2 + 2xy^5 + y^6$.
Trying~to generate
$G$
by~$3$,
$x+y+z$,
$f_3$,
and~$f_5$
now leads to linear system of seven equations in six unknowns which is easily seen to be~unsolvable.
\eop
\end{ex}

\end{document}